**Udo Pachner (1947-2002) – A "Hidden Champion" in Mathematics**

**(Dedicated to Udo Pachner on the occasion of his 70$^{th}$ birthday on April 20, 2017)**

**Peter Kleinschmidt, University of Passau**

## Introduction

This paper is dedicated to Professor Udo Pachner who was a very close friend and colleague of mine from 1968 up to his early death in 2002.

It occurred to me that I should write this reminiscence when I received an invitation from the University of Bochum to give a keynote lecture in 2014 on the occasion of the 85$^{th}$ birthday of our joint thesis advisor Günter Ewald (who died in 2015).

While I was preparing my lecture and doing some literature research, I was surprised to find a substantial number of papers which used Pachner's main results.

Most of these references use his main result about "bistellar transformations", now commonly called "Pachner moves". The main application areas are p.l. topology and theoretical physics, especially loop quantum gravity theory.

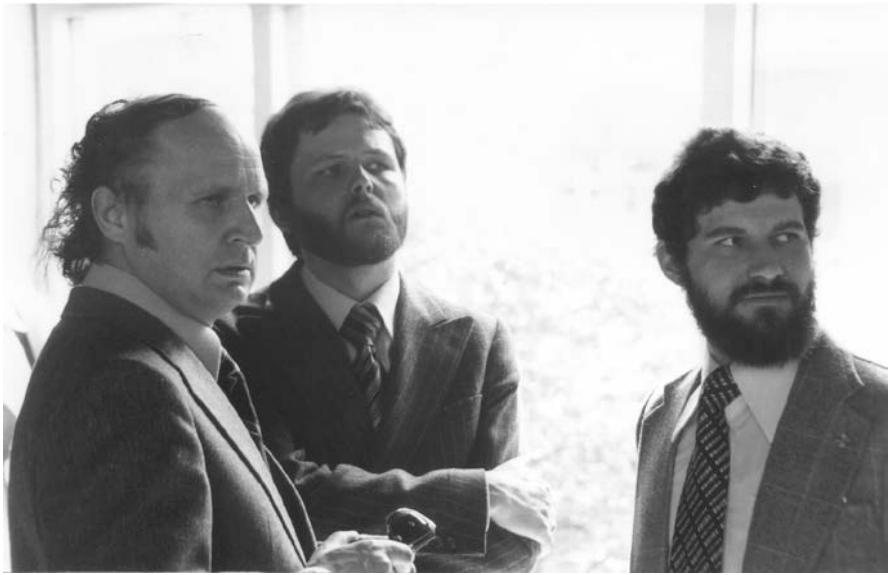

Günter Ewald, Peter Kleinschmidt, Udo Pachner (from left to right)
Colloquium on the occasion of the 90$^{th}$ birthday of Otto Haupt, Erlangen (1977)
Photo: Konrad Jacobs
Image Archive of the Mathematical Institute Oberwolfach, Nr. 1062

## The mathematics

Pachner moves are simple operations on simplicial complexes which can be described as follows:

Let *C* be a simplicial d-dimensional complex and A some k-dimensional simplex of *C* , 0≤ k ≤ d. The subcomplex Link(A,*C*) consists of all simplices L of *C* which are disjoint from A but for which A and L are faces of some other simplex of *C*. If in particular Link(A,*C*) is the boundary of some (d-k)-dimensional simplex B, B ∉ *C*, a Pachner move of dimension k on A consists of the following operation: Delete all simplices of *C* containing A and add the simplex B and all simplices which join B to faces of A. The new complex *C'* constructed this way is of course homeomorphic to *C*. *C'* admits the inverse Pachner move of dimension d-k on simplex B, thus creating *C* (simplicial isomorphisms are also always allowed).

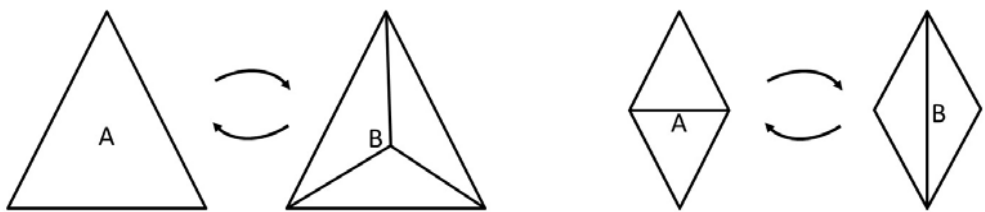

2-dimensional Pachner moves and their inverse

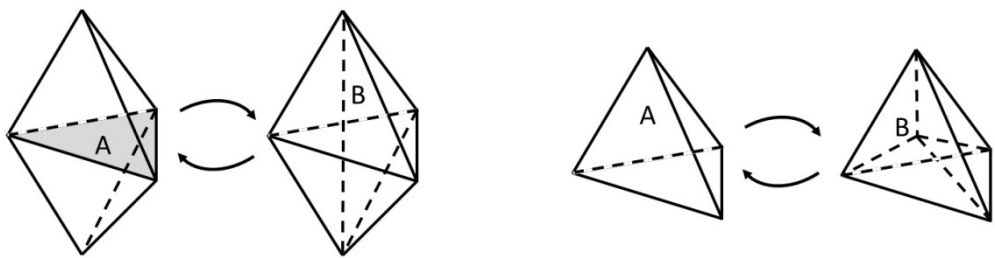

3-dimensional Pachner moves and their inverse

In [12], [13] und [14] Pachner proves that all closed simplicial manifolds which are homeomorphic can be transformed into each other by a finite sequence of Pachner moves (more precisely, he proved this for the category of p.l. manifolds). Ewald und Shephard had proved this result before for the boundary complexes of simplicial convex polytopes [6].

In the papers of Alexander and Newman [1], [10] and [11] it had already been proved that homeomorphic simplicial complexes can be transformed into each other by a finite sequence of "stellar" operations. Pachner proved a better result for a narrower class of objects (namely manifolds). Pachner moves can be interpreted as two subsequent stellar operations. The crucial property of Pachner moves is the fact that, in contrast to stellar operations, there are only d+1 types in dimension d. This makes it possible to identify invariants by proving the invariance only for a finite number of transformations in a fixed dimension. The reader may verify this easily by a very simple proof of the invariance of the Euler characteristic.

A first example of the use of this beneficial fact is the proof of the Turaev-Viro invariants in quantum theory, [8] and [15].

Motivated by this application, Lickorish honored Pachner's results in [8] and [9] by acknowledging their importance also in p.l. topology:

"In principle perhaps the basic theory of piecewise linear topology should be rewritten to accommodate Pachner's result, but market forces would inhibit such a project" [8].

Technically, the proofs of Pachner are quite complicated. They make heavy use of the results of Alexander. Pachner also introduces more transformations which he calls „elementary shellings". Using them he extends his results for closed manifolds to manifolds with boundary. Elementary shellings consist of additions or deletions of cells which induce Pachner moves on the boundary.

There are many purely mathematical papers which deal with Pachner moves, especially in p.l. topology. Lickorish's papers [8] and [9] provide very good insights.

A good example of its consequences is the following: handle structures which can be associated to different triangulations are connected by a sequence of additions or deletions of pairs of handles. This theory has its main importance in differential geometry, especially in Smale's h-cobordism theory.

Nevertheless it is helpful in p.l. topology that the existence of handle structures is obvious. This allows a more intuitive access, especially in low dimensions.

## The applications

The most quotations of Pachner moves can be found in loop quantum gravity theory. Here they have almost become a standard tool.

I am not a physicist, but I made some efforts to find relevant sources which provide a good point of entry for the non-expert and many further references for the interested expert. This is why I will restrict myself to some short explanations. The expert may forgive me for my choice of subject matter.

Loop quantum gravity theory, similar to string theory, is an attempt to find a unifying framework of quantum physics and general relativity theory.

Its basic assumption is that space or spacetime is not continuous but is represented by discrete entities in the order of magnitude of Planck length or Planck time, respectively.

Several approaches to the discretization have been established. One of them is a simplicial decomposition of space by which e.g. Pachner moves can be used to implement the dynamics of the systems.

The first applications of Pachner moves seem to have been the above-mentioned proof of the Turaev-Viro invariants in quantum theory. They correspond to certain state sums for arbitrary triangulations of a 3-manifold [15]. State sums are functions which serve as a tool for deriving physical state variables in quantum theory, as well as in other areas of physics like thermodynamics. In loop quantum gravity theory they are defined in a purely combinatorial way by assigning weights to the cells [15].

In [3] Pachner moves have been used to generalize the Dijkgraaf-Witten invariants [5] to dimension 4.

Further invariants for which Pachner moves proved essential can be found e.g. in [2], [4] and [7].

## Udo Pachner in private

Udo Pachner was a very amiable man who was much liked and highly esteemed both in his private and professional environments.

It was characteristic of him that he worked towards the proof of his main result in the most meticulous way. His wife recently told me that he used to say: "I always knew that it would work and that I would manage to find it".

It was very fortunate for me that I was able to witness the creation of his main work from close quarters. During the two years that he served as an acting professor at the University of Passau, he stayed at our house. We enjoyed discussing the complicated details of his proofs together.

His total dedication to his research work so impressed his wife and his daughter that they decided to install a plate at his gravestone which shows a Pachner move from his habilitation [12] thesis.

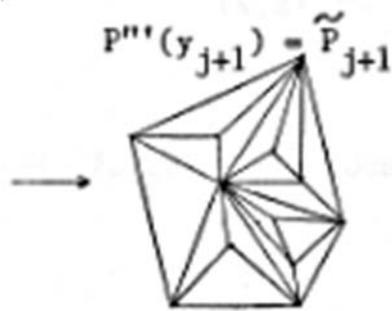 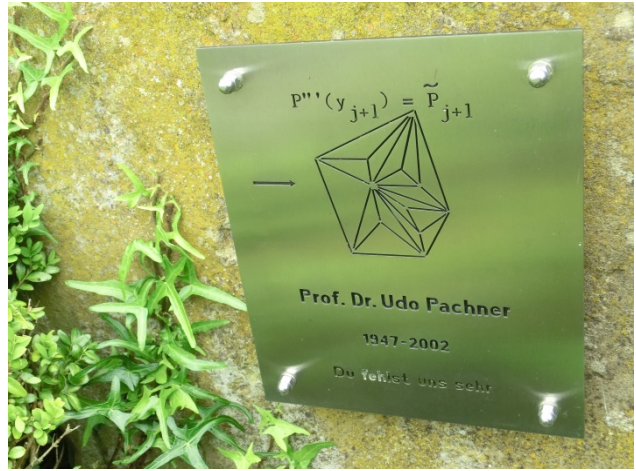

p. 25 of the habilitation thesis[12]     Gravestone, cemetery of Ahaus, Germany

Focusing so intensely on his work for his main result may have been a reason why he did not have the academic career he would have been capable of. During his last years he was in charge of a teaching institute at the Technical Academy in Ahaus, Germany, and also served as an adjunct professor at the Department of Mathematics of the University of Bochum.

You may call it tragic that he was not given the time to become aware of the importance of his results in loop quantum theory. However, he did read the above-mentioned positive remarks by Lickorish and was extremely pleased about them.

His wife and his daughter told me that despite his very ill health he died happy and content to be with his beloved family to the last. I am glad to have known him.

I am most grateful to Professor Lickorish for his most helpful remarks and to Udo's wife and daughter for providing some unpublished material.

Prof. Peter Kleinschmidt (Professor emeritus of the University of Passau)
Watzmannsdorfer Ring 27
D-94136 Thyrnau
pkleinschmidt@t-online.de